\documentclass{amsart}[12pt,A4]
\usepackage[margin=2.54cm]{geometry}
\usepackage[T1]{fontenc}
\usepackage{amsthm,amsmath,amssymb,color,graphicx,url}
\usepackage{float}
\usepackage{tikz}
\usepackage{graphicx}
\usepackage{mathtools}
\usepackage{makecell}
\usepackage{cleveref}
\usepackage{pstricks-add}
\usepackage[utf8]{inputenc}

\usetikzlibrary{calc}

\newtheorem{theorem}{Theorem}
\newtheorem{lemma}[theorem]{Lemma}
\newtheorem{corollary}[theorem]{Corollary}

\newtheorem{proposition}[theorem]{Proposition}
\theoremstyle{definition}

\newcommand{\ZZ}{\mathbb{Z}}

\newcommand{\Cov}{\hbox{\rm Cov}}

\newcommand{\Pos}{\hbox{Pos}}

\newcommand{\E}{{\rm E}}

\newcommand{\w}{{w}}

\newcommand{\M}{{\mathcal{M}}}
\newcommand{\F}{{\mathcal{F}}}
\newcommand{\B}{{\mathcal{B}}}
\newcommand{\C}{\mathcal{W}}

\newcommand{\po}{\mathcal{P}}

\usepackage[shadow, colorinlistoftodos,textsize=tiny, obeyFinal]{todonotes}
\setuptodonotes{fancyline, backgroundcolor=gray!10,bordercolor=gray}

\newcommand{\Mon}{{\rm Mon}}

\makeatletter
\numberwithin{equation}{section}
\numberwithin{figure}{section}
\numberwithin{theorem}{section}

\makeatother

\numberwithin{equation}{section}
\numberwithin{figure}{section}

\title{Sparse groups need not be semisparse}
\author{Isabel Hubard}
\address{Institute of Mathematics, National Autonomous University of Mexico (IM UNAM), 04510 Mexico City, Mexico}
\email{isahubard@im.unam.mx}
\author{Micael Toledo}
\address{}
\email{}

\begin{document}
\maketitle
\begin{abstract}

In 1999 Michael Hartley showed that any abstract polytope can be constructed as a double coset poset, by means of a C-group $\C$ and a subgroup $N \leq \C$.
Subgroups $N \leq \C$ that give rise to abstract polytopes through such construction are called {\em sparse}.
If, further, the stabilizer of a base flag of the poset is precisely $N$, then $N$ is said to be {\em semisparse}.
In \cite[Conjecture 5.2]{hartley1999more} Hartely conjectures that sparse groups are always semisparse.
In this paper, we show that this conjecture is in fact false: there exist sparse groups that are not semisparse. In particular, we show that such groups are always obtained from non-faithful maniplexes that give rise to polytopes.
Using this, we show that Hartely's conjecture holds for rank 3, but we construct examples to disprove the conjecture for all ranks $n\geq 4$.
\end{abstract}


\section{Introduction}


Abstract polytopes are partially ordered sets endowed with a rank function, which satisfy certain conditions that mimic the face poset of classical geometric polytopes (e.g. the Platonic Solids), tessellations (e.g. the cubical tessellation of 3-space), and some maps on surfaces (e.g. Klein's quartic). 
In this setting, {\em flags} correspond to the maximal chains of the poset and each flag $\Phi$ has a unique $i$-adjacent flag, which differs from $\Phi$ exactly on the element of rank $i$.
The flag graph of the polytope can be constructed as the edge-coloured graph whose vertices correspond to the flags and the edges of colour $i$ correspond to $i$-adjacent flags of the poset.

Maniplexes are edge-coloured graphs that simultaneously generalize maps on surfaces to structures of ranks higher than 3, and the flag graphs of abstract polytopes.
Every maniplex uniquely defines a partially ordered set; 
maniplexes whose vertices are in one-to-one correspondence with the maximal chains of their associated posets are called {\em faithful}. 
Conditions on when a faithful maniplex is the flag graph of a polytope are given in \cite{polimani}.
We will see (in \Cref{sec:ManiScher}) that maniplexes can be regarded as Schreier graphs and, thus, the graph theoretical conditions for polytopality can be translated to group theoretical conditions.

In 1999 Hartley (\cite{allpolytopes}) showed that every polytope is the quotient of a regular polytope.
This can be done by taking the Universal Polytope $\mathcal{U}$ of rank $n$ with automorphism group $\C^n$ together with some subgroup $N \leq \C$  and define
 the quotient $\po = \mathcal{U}/ N$ as the set of orbits of the faces of $\mathcal{U}$ under the action of $N$, where two such orbits are incident if they correspond to incident faces of $\mathcal{U}$. 
This is in fact equivalent to constructing $\po$ using a set of distinguished subgroups $\C_i \leq \C$, $i \in \{0,\ldots,n-1\}$ and construct 
the poset $\po=\po(\C,N)$ whose $i$-faces are the double cosets $\C_i wN$, and saying that $\C_iwN \leq \C_jvN$ if and only if $i\leq j$ and the two double cosets have a nonempty intersection. 
 Of course, not every $N\leq \C$ gives rise to a polytope. 
 Subgroups $N$ of $\C$ that produce polytopes in this way are called {\em sparse}.

Using the double coset construction,  Hartley showed that for each flag $\Psi$ of $\po$ there exists $w\in \C$ such that $\Psi$ consists of elements $\{\C_iwN \mid i\in \{0,\dots, n-1\}$, together with the minimal and maximal faces. 
Furthermore, setting $\Phi:=\{\C_iN \mid i \in\{0,1, \dots , n-1\} \} $, we may consider $N' \leq \C$ to be the stabilizer of $\Phi$ under the left action of $\C$.
It is straightforward to see that $N\leq N'$, but it is unclear if, even when $\po$ is in fact a polytope, $N$ and $N'$ coincide or not. 
Subgroups $N$ of $\C$ satisfying that $\po$ is a polytope and that $N=N'$ are called {\em semisparse}.

\if
As expected, one can express the conditions for a group $N\leq \C$ to be sparse and semisparse in terms of $N$ and several subgroups of $\C$. In fact, Hartley showed that the conditions to be semisparse are the following,
where $\C_{>i} = \langle r_j \mid j>i \rangle$ and $\C_{<i} = \langle r_j \mid j<i \rangle$:
\begin{enumerate}
    \item  $N^w \cap s_i\C_i = \emptyset$ for all $0\leq i \leq n-1$ and all conjugates $N^w$ of $N$ in $\C$.
    \item $\C_{>i} \cap \C_{<j}N^w \subseteq \C_k(\C_{>i}\cap \C_{<j}) N^w$, for all $w\in \C$ and $0\leq i<k<j \leq n-1$
\item $\C_{>i} \cap \C_{<i+1}N \subseteq N$, for all $0\leq i< n-1$.
\end{enumerate}
Moreover, if $N\leq \C$ satisfies the first two conditions above, then $N$ is in fact a sparse group.
\fi

Of course, if a group is semisparse, then it is sparse. However, the converse is an open question that we address in this paper: are sparse groups semisparse?
Hartley conjectures (see \cite[Conjecture 5.2]{hartley1999more}) that this is the case; after all, the idea behind his double coset construction is precisely that one uses the stabilizer of the base flag to construct a polytope. Furthermore, he shows that given a sparse group $N\leq \C$, and the group $N'\leq \C$ in such a way that $N'$ is the stabilizer of the base flag $\{\C_i N \mid i\in\{0,1,\dots, n-1\}\}$ of $\po(\C,N)$, then the polytopes $\po(\C,N)$ and $\po(\C,N')$ are isomorphic.

In this paper, we show that Hartley's conjecture is in fact false: there exist sparse groups that are not semisparse. 
In fact, we show that sparse groups that are not semisparse are always obtained from non-faithful maniplexes whose partial order is a polytope.
Thus, we settle the conjecture by constructing maniplexes with these properties for all ranks $n>3$. 
Moreover, we also show that
in rank $3$, sparse and semisparse groups indeed coincide, making Hartley's conjecture true for $n=3$.

The paper is structured as follows. In Section \ref{sec:basic}, we define basic notions related to abstarct polytopes and maniplexes; furthermore, in \Cref{sec:ManiScher} we show that maniplexes can be regarded as Schreier graphs. 
In Section \ref{sec:polimani}, we explore the connection between maniplexes, polytopes, and sparse and semisparse groups; in particular, we show that sparse groups that are not semisparse correspond to unfaithful maniplexes whose corresponding poset is a polytope. As a corollary, we show that, for rank $3$, sparse groups are also semisparse.
In Section \ref{sec:bigconstruction} we construct a rank $4$ unfaithful polytopal maniplex, giving rise to our first group that is sparse but not semisparse. 
In Section \ref{sec:extensions} we extend the results of Section \ref{sec:bigconstruction} to higher ranks by means of colour-coded extensions and show that for all ranks $n>3$, there exist examples of sparse groups that are not semisparse.
Finally, Section 6 has some concluding remarks.

\section{Basic notions}\label{sec:basic}
\subsection{Abstract polytopes}

Abstract polytopes are combinatorial structures that generalize (the face lattice of) classical polytopes and tesellations. 
While in the classical theory polytopes have spherical facets and vertex figures, abstract polytopes permit a wide range of possibilities. 
Although the main interest in studying abstract polytopes has been that of their symmetries, they are interesting objects on their own.

Formally, an {\em (abstract) polytope of rank $n$} $\po$ (or simply an $n$-polytope) is a partially ordered set, endowed with a rank function from $\po$ to the set $\{-1, 0, \dots, n\}$ satisfying the following properties:
\begin{enumerate}
    \item $\po$ has a (unique) minimum element $F_{-1}$ and a (unique) maximal element $F_n$.
    \item All maximal chains, called {\em flags}, have exactly $n+2$ elements, including $F_{-1}$ and $F_n$.
    \item $\po$ is strongly flag-connected in the sense that given two flags $\Phi$ and $\Psi$ of $\po$, there exists a sequence $\Phi=\Phi_0, \Phi_1, \dots, \Phi_k = \Psi$ such that $\Phi_i$ and $\Phi_{i+1}$ differ in exactly one element and $\Phi\cap\Psi \subseteq \Phi_i$, for all $i$. 
    \item $\po$ satisfies the so-called {\em diamond condition}, which says that if $F,G \in \po$ are such that $F<G$ and $rank(F) - rank(G) = 2$, then there are exactly two faces $H \in \po$ that satisfy $F<H<G$ and $rank(H)=rank(F)+1$.
\end{enumerate}

The elements of a polytope $\po$ are called {\em faces}, and if a face has rank $i$ we may refer to it as a $i$-face of $\po$.
If two flags differ in exactly one face, they are said to be {\em adjacent}, and if the face in which they differ has rank $i$, then they are $i$-adjacent.
Given a flag $\Phi$ and $i \in \{ 0, 1, \dots, n-1\}$, the diamond condition implies that there is exactly one flag $\Phi^i$ such that $\Phi$ and $\Phi^i$ are $i$-adjacent. The set of all the flags of $\po$ is denoted by $\F(\po)$.

A {\em section} of a poset $\po$ is a closed interval of the partial order. That is, given $F,G \in \po$ with $F<G$, we define the section $G/F$ as the set $\{H\in\po\mid F\leq H\leq G \}$.
Note that  $\po$ is strongly flag-connected, if and only if every section of $\po$ (including $\po$ itself) is  {\em flag-connected}, in the sense that given any two flags, there is a sequence of adjacent flags from one to the other (see \cite{arp} for details).

Given a polytope $\po$, its {\em flag graph} $G_\po$ is the edge-coloured graph whose vertices are the flags of $\po$, and where there is an edge of colour $i$ between two vertices if the corresponding flags are $i$-adjacent. 
The universal string Coxeter group of rank $n$, denoted by $\C^n$, is the group generated by an ordered set of $n$ involutions $(r_0,r_1,\ldots,r_{n-1})$ satisfying that whenever two of them are not consecutive, they commute. That is $\C^n = \langle r_0, \ldots, r_{n-1} \rangle$ with the only defining relations being $r_i^2=1$ and $(r_ir_j)^2=1$ whenever $|i-j|>1$.

The group $\C^n$ acts (on the left) on the flags of any $n$-polytope as follows: the generators $r_i$ of $\C^n$ act by the rule $r_i\Phi = \Phi^i$, for every flag $\Phi$ and every $i \in \{0, \dots, n-1\}$. 
By the connectivity of $\po$, the action of $\C^n$ is transitive on the flags. 
Whenever the rank of the polytopes is implicit, we denote $\C^n$ simply as $\C$.
\subsection{Maniplexes}
A maniplex is a connected $n$-valent, properly $n$-edge-coloured graph with colours $0,1, \dots, n-1$, such that the $2$-factors of colours $i$ and $j$ are squares, whenever $\mid i-j \mid >1$. An edge of colour $i$ of a maniplex is called an $i$-edge. 
The vertices of a maniplex are usually called {\em flags}, since, clearly,  flag graphs of polytopes are maniplexes.

Let $r_i$ be the permutation of the flags of a maniplex $\M$ that interchanges the endvertices of every $i$-edge. 
Clearly $r_i$ is an involution for each $i \in \{0, \dots , n-1\}$. The group generated by the permutations $r_0, r_1, \dots, r_{n-1}$ is called the {\em monodromy (or connection) group } of $\M$, and we denote it by $\Mon(\M)$. 
The image of a flag $\Psi$ under the permutation $r_i$ is called the $i$-adjecent flag to $\Psi$, and is often denoted by $\Psi^i$. Inductively, $\Psi^{i_0,i_1, \dots, i_k} = (\Psi^{i_0, i_1, \dots i_{k-1}})^{i_k}$.

Note that the monodromy group of a maniplex $\M$ acts in the natural way on the flags of $\M$. That is, for a flag $\Psi$ of $\M$, we denote by $r_i\Psi$ the image of $\Psi$ under the permutation $r_i$, and thus, for $w=r_{i_0}r_{i_1}\dots r_{i_k}\in \Mon(\M)$, the flag $w\Psi$ is defined inducively as $r_{i_0} (r_{i_1}\dots r_{i_k})\Psi$. 
Note further that $r_ir_j\Psi = r_i(r_j\Psi)=r_i(\Psi^j)=(\Psi^j)^i=\Psi^{j,i}$.
In particular, this means that $r_{i_0} (r_{i_1}\dots r_{i_k})\Psi = (r_{i_1}\dots r_{i_k}\Psi)^{i_0}$, that is, $r_j w\Psi = (w\Psi)^j$ and $wr_i\Psi =w\Psi^i$. 
So, in general, $r_j w\Psi \neq wr_j \Psi$.

By the connectivity of $\M$ we see that $\Mon(\M)$ acts transitively on the flags of $\M$. 
However, it is worth observing that the action on $\Mon(\M)$ on the flags of $\M$ need not be an action by automorphims.
It is straightforward to see that the monodromy group of an $n$-maniplex $\M$ is a quotient of $\C=\C^n$, and that $\C$ acts transitively on the flags of $\M$.
Therefore, we abuse the notation and use $r_0, \dots, r_{n-1}$ to denote the generators of both $\C^n$ and $\Mon(\M)$.

Given $I\subset\{0,1, \dots, n-1\}$, we let $\M_I$ denote the graph resulting from $\M$ by deleting the edges of colour $j$, for all $j\notin I$. 
In particular, for $i\in \{0,1, \dots, n-1\}$, $\M_i:=\M_{\{i\}}$ is the perfect matching of edges of colour $i$. For convenience, we denote by $\M_{\bar{i}}$ to the graph $\M_{\{0,1, \dots, n-1\} \setminus \{i\}}$.
Given $\Phi \in \M$, $\M_I(\Phi)$ denotes the connected component of $\M_I$ that contains $\Phi$.
Thus, $\M_I(\Phi)$ contains all the flags $\Psi$ of $\M$ that can be reached from $\Phi$ by walking only along edges of colours in $I$. 
That is, $\M_I(\Phi) = \{w\Phi \mid w\in \langle r_i \mid i\in I \rangle \}$.

Let $\M$ and $\widetilde{\M}$ be two $n$-maniplexes. 
We say that ${\M}$ is a quotient of $\widetilde{\M}$ if there exists a mapping $\pi: \widetilde{\M} \to {\M}$, called a {\em quotient projection}, such that for every $\widetilde{\Psi} \in \M$ and $i\in\{0,1\dots, n-1\}$ we have $\pi(\widetilde{\Psi}^i) = (\pi(\widetilde{\Psi}))^i$.
This immediately implies that if $\pi$ is a quotient between $n$-maniplexes, then $\pi(\w\widetilde{\Psi}) = \w(\pi(\widetilde{\Psi}))$, for every $\w\in\C$

\begin{lemma}
    Let $\pi:\widetilde{\M} \to {\M}$ be a quotient of $n$-maniplexes such that there exists $\Phi \in \M$ with $\pi^{-1}(\Phi)$ containing more than one element. 
    Then, $Stab_{\C}(\widetilde{\Phi}) < Stab_{\C}(\Phi)$, where $\widetilde{\Phi} \in \pi^{-1}(\Phi)$.
\end{lemma}

\begin{proof}
    To show that $Stab_{\C}(\widetilde{\Phi}) \leq Stab_{\C}(\Phi)$, let $\w\in\C$ be such that $\w\widetilde{\Phi} = \widetilde{\Phi}$.
    Since $\pi$ is a quotient, then 
    $\w\Phi = \w(\pi(\widetilde{\Phi}))=\pi(\w\widetilde{\Phi})=\pi(\widetilde{\Phi})=\Phi$,
    implying that $\w\in Stab_{\C}(\Phi)$.

    Now, let $\widetilde{\Psi} \in \pi^{-1}(\Phi)$ with $\widetilde{\Psi}\neq\widetilde{\Phi}$, and let $\w\in\C$ be such that $\w\widetilde{\Phi} = \widetilde{\Psi}$. 
    Then $\w\notin Stab_{\C}(\widetilde{\Phi})$.
    However, $\w\Phi = \w(\pi(\widetilde{\Phi})) =\pi(\w\widetilde{\Phi}) = \pi(\widetilde{\Psi}) = \Phi$, so we have $\w\in Stab_{\C}(\Phi)$.
\end{proof}

\subsection{Maniplexes and Schreier graphs}\label{sec:ManiScher}
Let $G$ be a group generated by involutions $r_0, r_1, \dots, r_{n-1}$ and let $H\leq G$ be a subgroup of $G$ such that $r_i, r_ir_j \notin H$, for all $i,j\in \{0, \dots, n-1\}$ with $i\neq j$. 
The {\em Schreier graph of $G$ with respect to $H$ and the set $R=\{r_0, r_1, \dots, r_{n-1} \}$} is the graph whose vertices are the cosets $gH$, with $g\in G$ and there is an edge between $g_1H$ and $g_2H$ if there exists $i\in \{0, \dots, n-1\}$ such that $g_2H = r_ig_1H$. The Schreier graph of $G$ with respect to $H$ and the generating set $R$ is denoted by $\Delta(G, H, R)$.

Note that $r_iH = r_jH$ if and only if $r_ir_j^{-1} = r_ir_j \in H$, which implies that $\Delta(G, H, R)$ is a well-defined simple graph of valency $n$. 
Furthermore, by colouring each edge $(gH, r_igH)$ with the colour $i$, we obtain an $n$-edge-coloured graph, which implies that $\Delta(G, H, R)$ is a maniplex whenever the $2$-factors of colours $i$ and $j$ are squares for $\mid i-j \mid >1$. 
This is, $\Delta(G, H, R)$ is a maniplex whenever $r_i, r_ir_j \notin H$, and $r_ir_j$ is an involution, if $\mid i-j \mid >1$.

\begin{lemma}\label{maniplexascosetgraph}
Let $\M$ be an $n$-maniplex, let $\C = \langle r_0, r_1, \dots, r_{n-1}\rangle$ be the universal Coxeter group of rank $n$, and let $\Phi$ be a base flag of $\M$. Consider $N:=Stab_{\C}(\Phi)$, and $\Delta:=\Delta(\C, N, \{r_0, r_1, \dots, r_{n-1}\})$. 
Then $\Delta \cong \M$ as coloured graphs.    
\end{lemma}

\begin{proof}

Start by recalling that since $\C$ is transitive on the flags of $\M$, then the set of flags of $\M$ can be written as $\{w\Phi \mid w\in \C\}$.
Let $\varphi:\Delta \to \M$ be such that $\varphi(wN) = w \Phi$. We will show that $\varphi$ is the desired isomorphism.

First, suppose that $wN=w'N$. Then we have $w^{-1}w'\in N$, which implies $w^{-1}w'\Phi = \Phi$ and thus $\varphi(wN)=w\Phi = w'\Phi = \varphi(w'N)$. Therefore, $\varphi$ is a well-defined function. 
By the transitivity of $\C$ on the flags of $\M$, we have that $\varphi$ is onto.

Next, let $w_1,w_2\in \C$ be such that $\psi(w_1N)=\psi(w_2N)$. 
Thus, $w_1\Phi = w_2\Phi$, so $w_2^{-1}w_1 \in Stab_{\C}(\Phi) = N$.
Therefore, $w_2N=w_1N$ and $\varphi$ is a bijection.

Finally, the vertices $wN$ and $r_iwN$, go to the flags $w\Phi$ and $r_iw\Phi =(w\Phi)^i$, respectively, implying that a $i$-edge of $\Delta$ indeed goes to an $i$-edge of $\M$.
    
\end{proof}

Note that not all groups $N\leq \C$ give rise to maniplexes. 
For example, if $r_ir_j \in N$, then between the vertices $N$ and $r_iN$ there would be two edges (the one corresponding to $r_i$ and the one corresponding to $r_j$).
The following result is straightforward.
\begin{lemma}\label{Nmani}
    Let $N\leq \C$. Then $\Delta(\C, N, \{r_0, r_1, \dots, r_{n-1}\})$ is a maniplex if and only if $r_i, r_ir_j \notin N$, for all $i,j \in \{0,1,\dots, n-1\}.$
\end{lemma}

\section{Some connections between maniplexes and polytopes}\label{sec:polimani}
%
Garza-Vargas and Hubard showed in \cite{polimani} that given a maniplex $\M$, one can construct a ranked poset $\Pos(\M)$, in a natural way, as follows.
The elements of rank $i$ of $\Pos(\M)$ are the connected components of $\M_{\bar{i}}$ and two elements are incident if their corresponding connected components have non-empty intersection. The poset $\Pos(\M)$ then encodes information regarding incidence between the distinct faces of $\M$. When $\Pos(\M)$ is a polytope, we say that $\M$ is {\em polytopal}.  In the particular case when $\M$ is the flag-graph of a polytope $\po$, $\M$ is clearly polytopal as $\Pos(\M)$ is isomorphic to $\po$.
Of course, not all maniplexes are polytopal (see \cite{polimani} for conditions on a maniplex to be polytopal).

Note that there exists a surjective {\em flag function} $f \colon \M \to \Pos(\M)$ mapping every flag of $\M$ to its corresponding maximal chain in $\Pos(\M)$. 
In fact, given $\Phi \in \M$, $\Phi$ belongs to the faces of $\Pos(\M)$ that, as connected components of $\M$, contain $\Phi$. In other words, $f(\Phi)$ is the maximal chain of $\Pos(\M)$ that contains the connected components $\M_{\bar{i}}(\Phi)$, for $i\in \{0,1, \dots, n-1\}$.

We say $\M$ is {\em faithful} if its flag function is  injective.
That is, the chain $\M_{\bar{i}}(\Phi)$ of $\Pos(\M)$ is identified only with the flag $\Phi$. In other words, $\cap_{i=0}^{n-1} \M_{\bar{i}}(\Phi) = \Phi$, for all $\Phi\in\M$.
Accordingly, we say that $\M$ is {\em unfaithful} if its flag function is not injective. 
In this case, there exists $\Phi, \Psi \in \M$, with $\Phi\neq\Psi$, such that $f(\Phi)=f(\Psi)$, which implies that $\Psi, \Phi \in \cap_{i=0}^{n-1} \M_{\bar{i}}(\Phi)$.

A set $X$ of flags of an unfaithful maniplex $\M$ is said to be unfaithful if for every $\Phi,\Psi \in X$, both $\Phi$ and $\Psi$ are mapped to the same maximal chain in $\Pos(\M)$ under $f$. 

Observe that, in principle, the faithfulness and the polytopality of a maniplex need not be related. 
In particular, there are examples of faithful maniplexes that are polytopal (all the flag graphs of polytopes), as well as faithful maniplexes that are not polytopal (for example, the flag graph of the toroidal map $\{4,4\}_{(1,1)}$).
There are also examples of unfaithful maniplexes that are not polytopal (for example, the flag graph of the toroidal map $\{4,4\}_{(1,0)}$).
In Sections \ref{sec:bigconstruction} and \ref{sec:extensions} we present examples of unfaithful maniplexes of rank $n>3$ that are polytopal. To the best of our knowledge, these are the first examples of this phenomenon. 
In this section, we show that, in fact, all rank $3$ polytopal maniplexes are faithful. 

\begin{lemma}\label{lemma:02pairs}
A $3$-maniplex $\M$ is unfaithful if and only if there are unfaithful pairs of the form $\{\Phi,\Phi^0\}$ and $\{\Psi,\Psi^2\}$
\end{lemma}

\begin{proof}
Trivially, if $\M$ has an unfaithful pair, then $\M$ is itself unfaithful. 
Suppose that $\M$ is unfaithful. That is, there exist flags $\Phi$ and $\Psi$ that lie in the same $i$-face for all $i \in \{0,1,2\}$. 
In particular, $\Phi$ and $\Psi$ lie in the same $1$-face. If $\Psi = \Phi^j$ with $j \in \{0,2\}$, then we are done. Thus, since a $1$-face of $\M$ is a $4$-cycle of alternating colours $0$ and $2$, we have $\Psi = \Phi^{02}$. 
Now, $\Phi^0$ lies on the same $j$-face than $\Phi$, for $j \in \{1,2\}$. 
Also, $\Phi^0 = \Psi^2$ and thus $\Phi^0$ lies on the same $0$-face as $\Psi$, but by hypothesis $\Psi$ and $\Phi$ lie in the same $0$-face. 
That is, $\Phi$ and $\Phi^0$ lie on the same $i$-face for $i \in \{0,1,2\}$. Therefore, $\{\Phi,\Phi^0\}$ is an unfaithful pair. A similar argument shows that $\{\Psi,\Psi^2\}$ is an unfaithful pair.  
\end{proof}

\begin{proposition}\label{unfaithfulrank3}
If $\M$ is an unfaithful $3$-maniplex, then $\Pos(\M)$ is not an abstract polytope.
\end{proposition}

\begin{proof}
Suppose $\Pos(\M)$ is an abstract polytope and let $f$ be the flag function of $\M$. Since $\M$ is unfaithful, there exists an unfaithful pair $\{\Phi,\Phi^2\}$. Let $C= f(\Phi) = f(\Phi^2)$ and let $C'$ be the maximal chain in $\Pos(\M)$ that differs from $C$ only in their $2$-face. Such a pair of chains is guaranteed to exist because $\Pos(\M)$ is an abstract polytope. Let $\Psi \in f^{-1}(C')$ and observe that $\Psi$ must be in the same $1$-face as $\Phi$ and $\Phi^2$. But then $\Psi$ has no other choice but to be one of $\Phi^0$ or $\Phi^{2,0}$ and thus $\Psi$ is in the same $2$-face as either $\Phi$ or $\Phi^2$, contradicting that $C$ and $C'$ differ on their $2$-face. We conclude that $\Pos(\M)$ is not an abstract polytope.
\end{proof}

Note that to show \Cref{lemma:02pairs} and \Cref{unfaithfulrank3} we strongly use the fact that a $1$-face of a $3$-maniplex has exactly $4$ flags. This is no longer the case for $n$-maniplexes with $n>3$. In fact, the result fails for $n>3$, as we show in \Cref{sec:bigconstruction}.

\subsection{Sparse groups and maniplexes}

As mentioned in the introduction, in \cite{allpolytopes} Hartley constructs posets by means of double cosets. 
To do so, he uses the Universal Coxeter group $\C=\langle r_0, r_1, \dots, r_{n-1} \mid r_i^2, (r_ir_j)^2=1\ \mathrm{whenever} \ |i-j|>1 \rangle$, and for each $N\leq \C$ constructs
the poset $\po(\C,N)$ whose $i$-faces are the double cosets $\C_iwN $, and says that $\C_iwN \leq \C_jvN$ if and only if $i\leq j$ and the two double cosets have a nonempty intersection, where $\C_i = \langle r_j \mid j\neq i\rangle$, for each $i \in \{0, 1, \dots, n-1\}$.
Although $\po(\C,N)$ need not be a polytope, it is always a ranked poset with maximal chains all of the same length. 

We say that a group $N\leq \C$ is {\em sparse} if $\po(\C,N)$ is a polytope. Even though this is not the way sparse groups are defined in \cite{mcmullen1994quotients}, the results of \cite{allpolytopes,hartley1999more} imply the equivalence in the definitions (see also \cite{hartley2006simpler}).

In this section we  see connections between the poset $\po(\C,N)$ and maniplexes.
Given $N\leq \C$ such that $r_i, r_ir_j \notin N$, for all $i,j\in \{0,1,\dots, n-1\}$, we have seen how to construct the maniplex $\M=\Delta(\C,N,\{r_0,r_1, \dots, r_{n-1}\})$, and thus we can construct the poset $\Pos(\M)$.
On the other hand, we can construct the
 poset $\po=\po(\C,N)$, and thus its flag graph $\M(\po)$.

\begin{proposition}\label{equalposets}
Let $N\leq \C$ be such that $r_i, r_ir_j \notin N$ for all $i,j \in\{0,1,\dots, n-1\}$.
Then $\Pos(\M)=\po(\C,N)$, where $\M=\Delta(\C, N , \{r_0, r_1, \dots, r_{n-1} \})$
\end{proposition}

\begin{proof}
First note that by \Cref{Nmani} $\M$ is a maniplex, and hence $\Pos(\M)$ is a well-defined poset.
Observe also that the flags of $\M$ are of the form $wN$, with $w\in\C$.

Recall that the $i$-faces of $\Pos(\M)$ are the connected components of $\M_{\bar{i}}$.
To show that $\Pos(\M)=\po(\C,N)$ we first need to show that the flags in connected component of a flag $wN$ of $\M_{\bar{i}}$ is precisely the set of flags $\C_i wN$.

Two flags $wN$ and $uN$ of $\M$ are the same connected component of $\M_{\bar{i}}$ if and only if there is a path between $wN$ and $uN$ without edges of colour $i$.
That is, $wN$ and $uN$ are in the same connected component of $\M_{\bar{i}}$ if and only if there exists $p \in \C_i$ such that $pwN=uN$.
Therefore, the flags in $\{pwN \mid p \in \C_i\}$ are precisely the flags in the connected component of $\M_{\bar{i}}$ containing $wN$.
That is, the connected component of $\M_{\bar{i}}$ containing $wN$ is precisely $\C_i wN$, and thus the $i$-faces of $\Pos(\M)$ are exactly the same as the $i$-faces of $\po=\po(\C,N)$.

Now let $\C_iwN \leq_\po \C_jvN$.
By \cite[Lemma 2.2]{allpolytopes}, there exists $u \in \C$ such that $\C_iwN = \C_i u N$ and $\C_jvN = \C_j u N$. This implies that $uN$ is in the connected component of $\M_{\bar{i}}$ containing $wN$ as well as in the connected component of $\M_{\bar{j}}$ containing $vN$. Thus those two connected components have non empty intersection. Therefore $\C_iwN \leq_\M \C_jvN$.

Conversely, if the connected component of $\M_{\bar{i}}$ containing $wN$ and the connected component of $\M_{\bar{j}}$ containing $vN$ have non empty intersection, then there exists a flag $uN$ in both of them. 
Thus, $uN \in \C_i wN \cap \C_j vN$, implying 
$\C_iwN \leq_\po \C_jvN$.

This means that the sets of $\Pos(\M)$ and $\po(\C,N)$ are equal and that two elements of one are incident if and only if they are incident in the other one. The theorem follows.
\end{proof}

Recall, on one hand, that a maniplex $\M$ is polytopal if and only if $\Pos(\M)$ is a polytope; and on the other hand, that $N\leq \C$ is sparse if and only if $\po(\C,N)$ is a polytope. The following result is a direct consequence of these two definitions, \Cref{Nmani} and \Cref{equalposets}.

\begin{corollary}\label{coro:sparse-polytopal}
Let $N\leq \C$.
Then $N$ is a sparse group if and only if $\Delta(\C, N , \{r_0, r_1, \dots, r_{n-1} \})$ is a polytopal maniplex.
\end{corollary}

\subsection{Sparse groups that are not semisparse.}

Let $N\leq \C$ be a sparse group, and let $\po=\po(\C,N)$ its corresponding polytope with base flag $\Phi = \{\C_iN \mid i\in \{0,1,\dots, n-1 \}$.
Let $N' = Stab_{\C}(\Phi)$.
We say that the group $N$ is {\em semisparse} if $N=N'$.
The following results imply that to find sparse groups that are not semisparse, we need to find unfaithful maniplexes that are polytopal. 

\begin{proposition}\label{notsemisparseimplesunfaithful}
    Let $N\leq \C$ be a sparse group that is not semisparse. Then the maniplex $\M=\Delta(\C,N,\{r_0,r_1, \dots, r_{n-1}\})$ is unfaithful and polytopal.
\end{proposition}

\begin{proof}
    By \Cref{coro:sparse-polytopal} $\M$ is a polytopal maniplex. In particular this implies that $r_i, r_ir_j \notin N$. Furthermore, this also implies that
 we only need to show that $\M$ is unfaithful.

    Since $N$ is sparse, then the poset $\po(\C,N)$ is a polytope and, if $\Phi=\{\C_iN\mid i\{0,1,\dots, n-1\}\}$ and $N'=Stab_{\C}(\Phi)$, the fact that $N$ is not semisparse implies that $N$ is a proper subgroup of $N'$.
    As pointed out before, in \cite{hartley1999more} Hartley shows that $\C_iwN = \C_iwN'$ for all $w\in \C$ and $i\in\{0,1,\dots, n-1\}$, implying that  $\po(\C,N) =\po(\C,N')$, which in turn implies that $N$ is sparse if and only if $N'$ is sparse.
    Thus, $r_i, r_ir_j \notin N'$.

    Consider $\M':=\Delta(\C,N',\{r_0,r_1, \dots, r_{n-1}\})$. 
    By \Cref{equalposets} and the above observation we have that $\Pos(\M) = \po(\C,N) = \po(\C,N') = \Pos(\M')$.

    Now define $\vartheta: \M \to \M'$ in such a way that $\vartheta(wN) = wN'$.
    Then $\vartheta$ is well-defined, as if $w_1N = w_2N$, then $w_2^{-1}w_1\in N < N'$, implying that $w_1N' = w_2N'$, and thus $\vartheta(w_1N)=\vartheta(w_2N)$. 
    Moreover, $\vartheta$ is clearly onto.

    Note further that $(wN)^i=r_iwN$, which implies that $\vartheta((wN)^i) = \vartheta(r_iwN) = r_iwN' = r_i \vartheta(wN) = (\vartheta(wN))^i$, implying that $\vartheta$ is a quotient between the maniplexes $\M$ and $\M'$.

    Now, since $N$ is a proper subgroup of $N'$, there exists $w_0\in N'\setminus N$. Then $w_0N\neq N$, but $\vartheta(w_0N)=w_0N'=\vartheta(N)$, which implies that $\vartheta$ is not one-to-one. 
    
    Thus, $\M$ has more flags than $\M'$, while $\Pos(\M)=\Pos(\M')$. This means that $\M$ is unfaithful (in fact, the flags $N$ and $w_0N$ are an unfaithful pair).
\end{proof}


\begin{lemma}\label{propersubgroup}
    Let $\widetilde{\M}$ be
  a 
    maniplex such that $\po:=\Pos(\widetilde{\M})$ is an abstract polytope, and let $\M:=\M(\po)$ be the flag graph of $\po$. Let $\widetilde{\Psi}$ be a flag of $\widetilde{\M}$. Then:
    \begin{enumerate}
        \item $\M$ is a quotient of $\widetilde{\M}$,
        \item $Stab_{\C}(\widetilde{\Psi}) \leq Stab_{\C}(\Psi)$, where $\widetilde{\Psi}$ projects to the flag $\Psi$ of $\M$; furthermore, if $\widetilde{\M}$ is unfaithful and $\widetilde{\Psi}$ belong to an unfaithful set, then $Stab_{\C}(\widetilde{\Psi}) < Stab_{\C}(\Psi)$
        \item $Stab_{\C}(\Psi) = Stab_{\C}(f(\widetilde{\Psi}))$, where  $f: \widetilde{\M} \to \Pos(\widetilde{\M})$ is the flag function.
    \end{enumerate}
\end{lemma}

\begin{proof}
Note that by using the flag function from $\widetilde{\M}$ to $\po$ and the bijection between the flags of $\po$ and its flag graph we can obtain a quotient map $\vartheta:\widetilde{\M} \to \M$. 
Furthermore, $\vartheta$ is onto, since every flag of $\M$ corresponds to a flag of $\po$, which in turn corresponds to at least one flag of $\widetilde{\M}$. This settles item (1).

Let $\widetilde{\Psi}\in \widetilde{\M}$ and  $\vartheta(\widetilde{\Psi})=\Psi$.
If $w\in Stab_{\C}(\widetilde{\Psi})$, then $w(\widetilde{\Psi}) = (\widetilde{\Psi})$, which implies that $\Phi = \vartheta(\widetilde{\Psi})=\vartheta(w(\widetilde{\Psi}))= w\vartheta((\widetilde{\Psi}))=w\Phi$. Hence $Stab_{\C}(\widetilde{\Psi})\leq Stab_{\C}(\Phi)$.  
Now, if $\widetilde{\M}$ is unfaithful, and $\{\widetilde{\Psi}, \widetilde{\Phi}\}$ is an unfaithful set, then $\vartheta(\widetilde{\Psi})=\vartheta(\widetilde{\Phi})$.
Let $u \in \C$ be such that $u\widetilde{\Psi}=\widetilde{\Phi}$. Then $u\notin Stab(\widetilde{\Psi})$, but $\Phi = \vartheta(\widetilde{\Phi}) = \vartheta(u\widetilde{\Psi})=u\vartheta(\widetilde{\Psi})=u\Phi$, so $u \in Stab_{\C}(\Phi)$, settling item (2).

Item (3) is straightforward, as the flags of $\M(\po)$ correspond to the flags of $\po$ and the action of $\C$ is the same in these two sets. 
\end{proof}

\begin{proposition}
\label{unfaithfulpolytopal}
    Let $\widetilde{\M}$ be an unfaithful maniplex that is polytopal, and let $\widetilde{\Psi} \in \widetilde{\M}$. Then $N=Stab_{\C}(\widetilde{\Psi})$ is a sparse group that is not semisparse.   
\end{proposition}

\begin{proof}
    Since $\widetilde{\M}$ is polytopal, by \Cref{coro:sparse-polytopal}, $N$ is a sparse group. Furthermore,
by \Cref{equalposets} we have $\Pos(\widetilde{\M}) = \po(\C,N)$.

Now $N$ is semisparse if and only if $N=Stab_{\C}(\Phi)$, where $\Phi$ is the flag $\{\C_iN \mid i\in \{0,1, \dots, n-1\}\}$.
By \Cref{propersubgroup}, $N < Stab_\C (f(\widetilde{\Psi}))$, where $f: \widetilde{\M} \to \Pos(\widetilde{\M})$ is the flag function. Thus, to show that $N$ is not semisparse, we need to show that $\Phi = f(\widetilde{\Psi})$. 

Since $\widetilde{\M} = \Delta(\C, N, R)$, the flag $\widetilde{\Psi}$ can be relabelled as the vertex $N$, and thus $f(N)$ is the maximal chain of $\Pos(\widetilde{\M})$ that contains all connected components $\widetilde{\M}_{\bar{i}}(N)$, for $i\in \{0,1,\dots, n-1\}$.
As pointed out in the proof of \Cref{equalposets}, for each $i$, the connected component of $\widetilde{\M}_{\bar{i}}$ that contains the flag $wN$ is precisely the set $\C_iwN$, so, in particular, this means that the connected component of $\widetilde{\M}_{\bar{i}}$ that contains the flag $N$ is precisely the set $\C_iN$. 
This implies that $\Phi = f(\widetilde{\Psi})$, and thus $N$ is not semisparse.
\end{proof}

\Cref{notsemisparseimplesunfaithful} and \Cref{unfaithfulpolytopal} imply the following theorem.

\begin{theorem}\label{teo:theequivalnce}
    Let $N \leq \C$ be a sparse group. Then $N$ is semisparse  if and only if the polytopal maniplex $\Delta(\C, N, \{r_0, r_1, \dots, r_{n-1}\})$ is faithful.
\end{theorem}

The above theorem, together with \Cref{unfaithfulrank3} imply the following result.

\begin{corollary}
    If $\C$ is the universal Coxeter group of type $[\infty,\infty]$ of rank $3$ and $N\leq \C$ is a sparse group, then $N$ is semisparse.
\end{corollary}

In the following sections we will see that the above result only holds for rank $n=3$. This is done by constructing polytopal maniplexes that are unfaithful for all $n>3$.

\section{An unfaithful polytopal $4$-maniplex}\label{sec:bigconstruction}
We have seen that in order to find a group that is sparse but not semisparse, we need to have an unfaithful maniplex that is polytopal.
In this section, we construct such a maniplex for $n=4$. 
In \Cref{sec:extensions} we then use this maniplex to construct unfaithful maniplexes of rank $n>4$ that are polytopal.

\subsection{Double covers from voltage assignments}\label{sec:construction}

Let $\Gamma$ be a graph with vertex set $V$ and edge set $E$. A {\em $\ZZ_2$-voltage assignment} is a function $\zeta\colon E \to \ZZ_2$. The {\em double cover} of $\Gamma$ relative to $\zeta$, denoted $\Cov(\Gamma,\zeta)$, is the graph with vertex set $V \times \ZZ_2$ and where two vertices $(x,i)$ and $(y,j)$ are adjacent if and only if $e:=\{x,y\}$ is an edge of $\Gamma$ and $j = i + \zeta(e)$.

Note that the mapping $\pi\colon \Cov(\Gamma,\zeta) \to \Gamma$ sending every vertex $(x,i)$ to $x$ is a quotient projection and thus $\Gamma$ is a quotient of $\Cov(\Gamma,\zeta)$. 
When $\Gamma$ is a rank $n$ maniplex, we can give a proper colouring to $\Cov(\Gamma,\zeta)$ by colouring every edge $e$ with the same colour as $\pi(e)$.
Note that this is not enough to guarantee that $\Cov(\Gamma,\zeta)$ is an $n$-maniplex: we need to ensure that $\Cov(\Gamma,\zeta)$ is connected and that the subgraphs induced by non-consecutive colours are unions of squares. 
The lemma \ref{lem:voltcon} is folklore, while Lemma \ref{lem:voltmani} follows from the definitions of a maniplex and a double cover.

\begin{lemma}
\label{lem:voltcon}
Let $\Gamma$ be a connected graph, let $\zeta$ be a $\ZZ_2$-voltage assignment, and let $\Gamma'$ be a connected subgraph of $\Gamma$. If the set of edges of $\Gamma'$ with non-trivial voltage is not a cut-set of $\Gamma'$, then $\pi^{-1}(\Gamma')$ is connected. 
\end{lemma}

\begin{lemma}
\label{lem:voltmani}
Let $\M$ be a maniplex (of rank $n>2$) and let $\zeta$ be a $\ZZ_2$-voltage assignment. 
Then $\Cov(\M,\zeta)$ is a maniplex if and only if the following hold:
\begin{enumerate}
\item the set of edges of $\M$ with non-trivial voltage is not a cut-set;
\item every $4$-cycle of $\M$ of alternating colours $i$ and $j$, with $|i-j|>1$, has an even number of edges with non-trivial voltage.
\end{enumerate}
\end{lemma}

\subsection{An unfaithful polytopal $4$-maniplex}

The goal of this construction is to obtain an unfaithful $4$-maniplex whose poset is an abstract polytope. 
To do so, we start with a faithful $4$-maniplex $\B$ that is polytopal and reflexible (that is, maximally symmetric).
Then, we construct a double cover $\B^*$ of $\B$ such that for all $i \in \{0,1,2,3\}$ and every $i$-face $F$ of $\B$, the lift $\pi^{-1}(F)$ is connected, and thus an $i$-face of $\B^*$.
Therefore, there is a bijection between the $i$-faces of $\B$ and the $i$-faces of $\B^*$ that preserves the order, which implies that $\B^*$ and $\B$ have isomorphic posets.

\begin{figure}[h!]
\begin{center}
\includegraphics[width=0.85\textwidth]{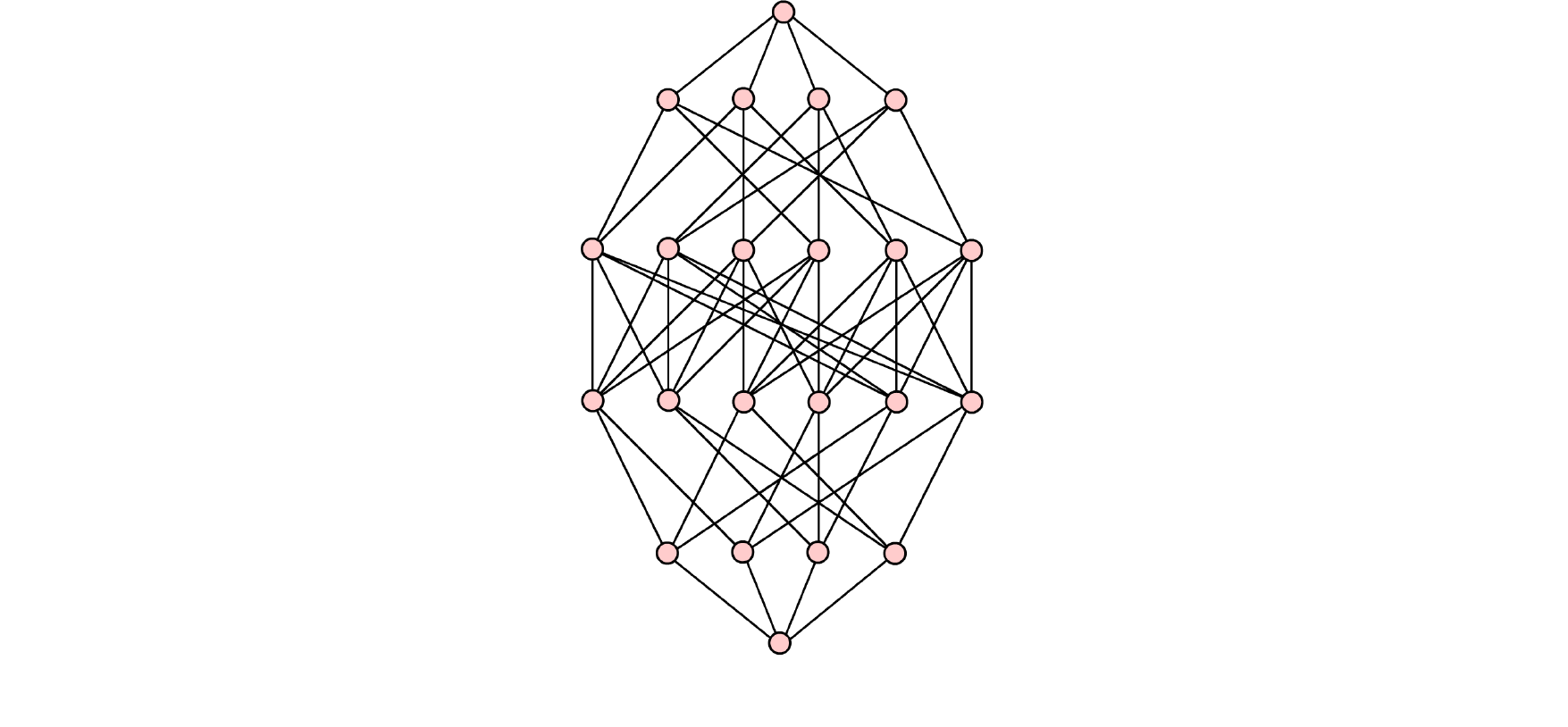}
\caption{The self-dual regular $4$-polytope $\po$.}
\label{fig:posP}
\end{center}
\end{figure}

Let $\po$ be the regular self-dual $4$-polytope $\{\{4,3\}_3,\{3,4\}_3\}$. This polytope is denoted by $\{4,3,4\}^*96$ in \cite{hartley2006atlas}. It is the universal polytope whose facets are hemi-cubes $\{4,3\}_3$, and whose vertex-figures are hemi-octahedra $\{3,4\}_3$, and has a total of $96$ flags. 
It has $4$ vertices, $6$ edges, $6$ polygons (or $2$ faces, all squares) and $4$ facets (all hemicubes), from which we can see that it is a flat polytope (that is, every vertex belongs to every facet).
The Hasse diagram of $\po$ can be seen in \Cref{fig:posP}.

Let $\B$ be the flag graph of $\po$. 
Then $\B$ is a faithful polytopal $4$-maniplex, and the flags of $\po$ are in one-to-one correspondence with the flags of $\B$.
To construct the double cover of $\B^*$, we need to define a $\ZZ_2$-voltage assignment for the edges of $\B$.
To do so, we use the Hasse diagram of $\po$ and define some sets of flags and $i$-edges of $\B$.
%

\begin{figure}[h!]
\begin{center}
\includegraphics[width=0.85\textwidth]{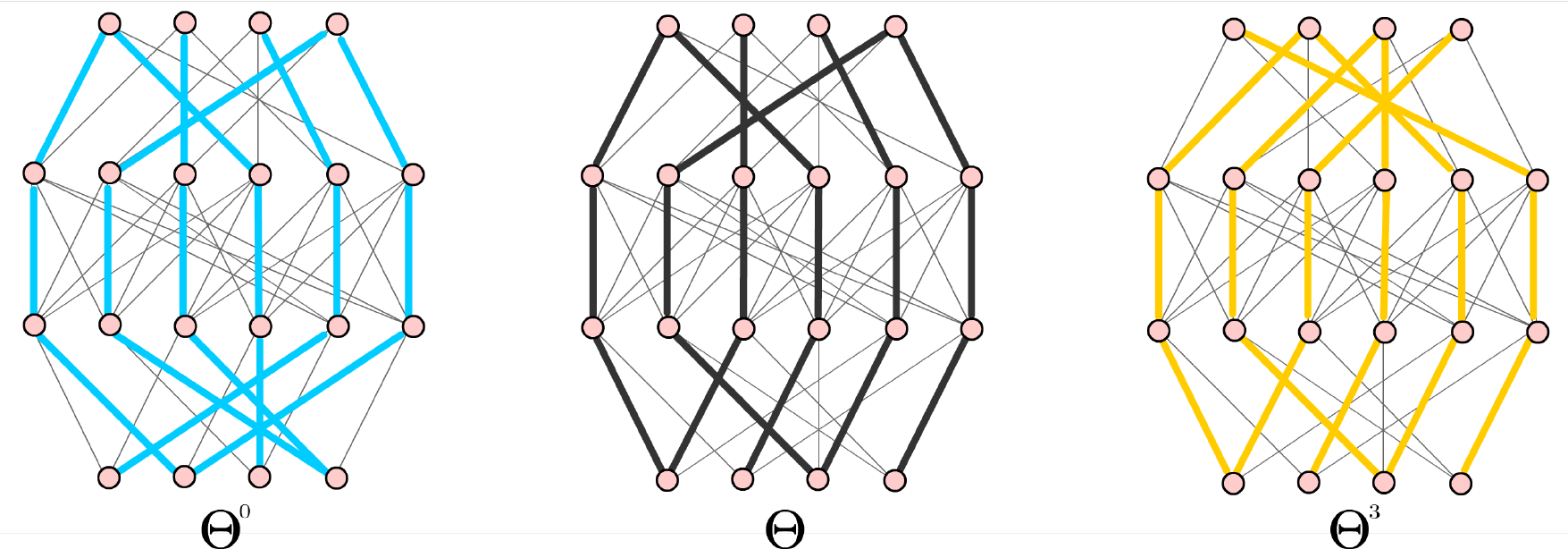}
\caption{The distinguished sets of flags $\Theta^0$, $\Theta$ and $\Theta^3$ (the minimal and maximal faces are not shown)}
\label{fig:posets}
\end{center}
\end{figure}

Let $\Theta$ be the set of flags of $\B$ corresponding to the maximal chains of $\po$ drawn in bold black edges, in the centre of Figure \ref{fig:posets}. For $i \in \{0,1,2,3\}$, let $\Theta^i = \{\Phi^i \mid \Phi \in \Theta \}$. The sets $\Theta^0$ and $\Theta^3$ are depicted (as chains on the poset $\po$) in the left-hand and right-hand side of Figure \ref{fig:posets}, respectively. 

We have chosen the flags of $\Theta$ so that $\Theta$ has certain desirable properties. 
In particular if $F$ is an $i$-face of $\B$ (that is, a connected component of $\B_{\hat{i}}$)
then the following hold:
\begin{enumerate}
\item[(A.1)] If $i \in \{1,2\}$, then $F$ contains exactly one flag $\Phi \in \Theta$;
\item[(A.2)] If $i \in \{0,3\}$, then $F$ contains either exactly one or exactly two flags in $\Theta$;
\item[(A.3)] If $i \in \{0,3\}$ and $F$ contains two flags in $\Theta$, then these flags lie on different $j$-faces for all $j \neq i$. Moreover, $F$ contains exactly one flag in $\Theta^i$;
\item[(A.4)] If $i \in \{0,3\}$ and $F$ contains only one flag in $\Theta$, then it contains exactly two flags in $\Theta^i$.
\end{enumerate}

\begin{figure}[h!]
\begin{center}
\includegraphics[width=0.85\textwidth]{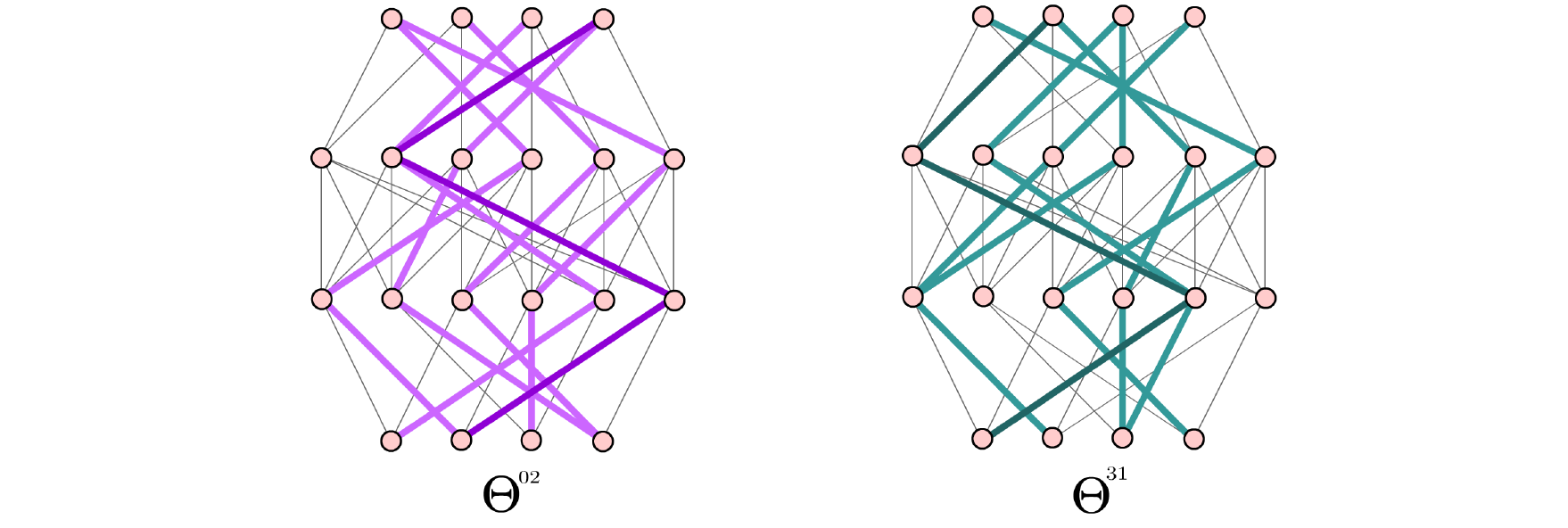}
\caption{The distinguished sets of flags $\Theta^{02}$ and $\Theta^{31}$ (the minimal and maximal faces are not shown)}
\label{fig:pos}
\end{center}
\end{figure}

To define a voltage assignment for $\B$, we define a set of edges, $E_\Theta$, to carry non-trivial voltages, as follows. 
First, for each $\Phi \in \Theta$, let $E_\Phi$ be the set of edges
$\{(\Phi,\Phi^0),(\Phi^0,\Phi^{02}),(\Phi,\Phi^3),(\Phi^3,\Phi^{31})\}$ (see \Cref{fig:pos} for the sets containing $\Phi^{02}$ and $\Phi^{13}$, for $\Phi \in \Theta$).

\begin{figure}[H]
\begin{center}
\includegraphics[width=0.35\textwidth]{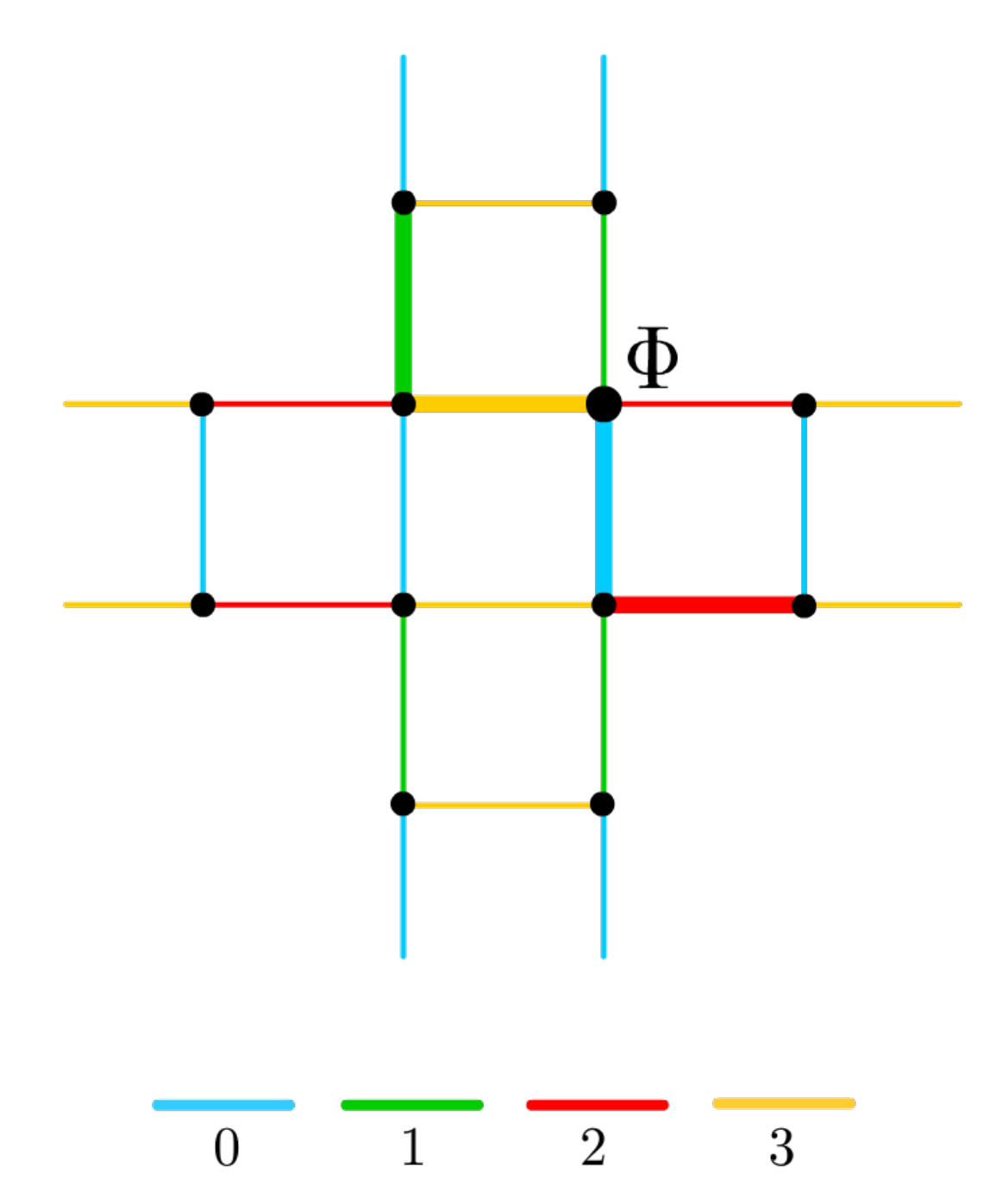}
\caption{Part of the intersection of a $1$- with a $2$-face of $\B$. A set $E_\Phi$ shown in thick edges.}
\label{fig:voltaje}
\end{center}
\end{figure}

Note that the edges of $E_\Phi$ are the edges of a $1302$-path that starts at $\Phi^{31}$ and ends at $\Phi^{02}$ (see Figure \ref{fig:voltaje}). 
Let $E_\Theta = \bigcup_{\Phi \in \Theta} E_\Phi$. 

We are now ready to define a $\ZZ_2$-voltage assignment $\zeta_\Theta: \E(\B) \to \ZZ_2$. We let

\[
\zeta_\Theta(e)=
 \begin{cases}
1 & \text{ if $e \in E_\Theta$}   \\
0 & \text{ otherwise.}
\end{cases}
\]

\begin{figure}[h]
    \centering
    \includegraphics[width=0.6\textwidth]{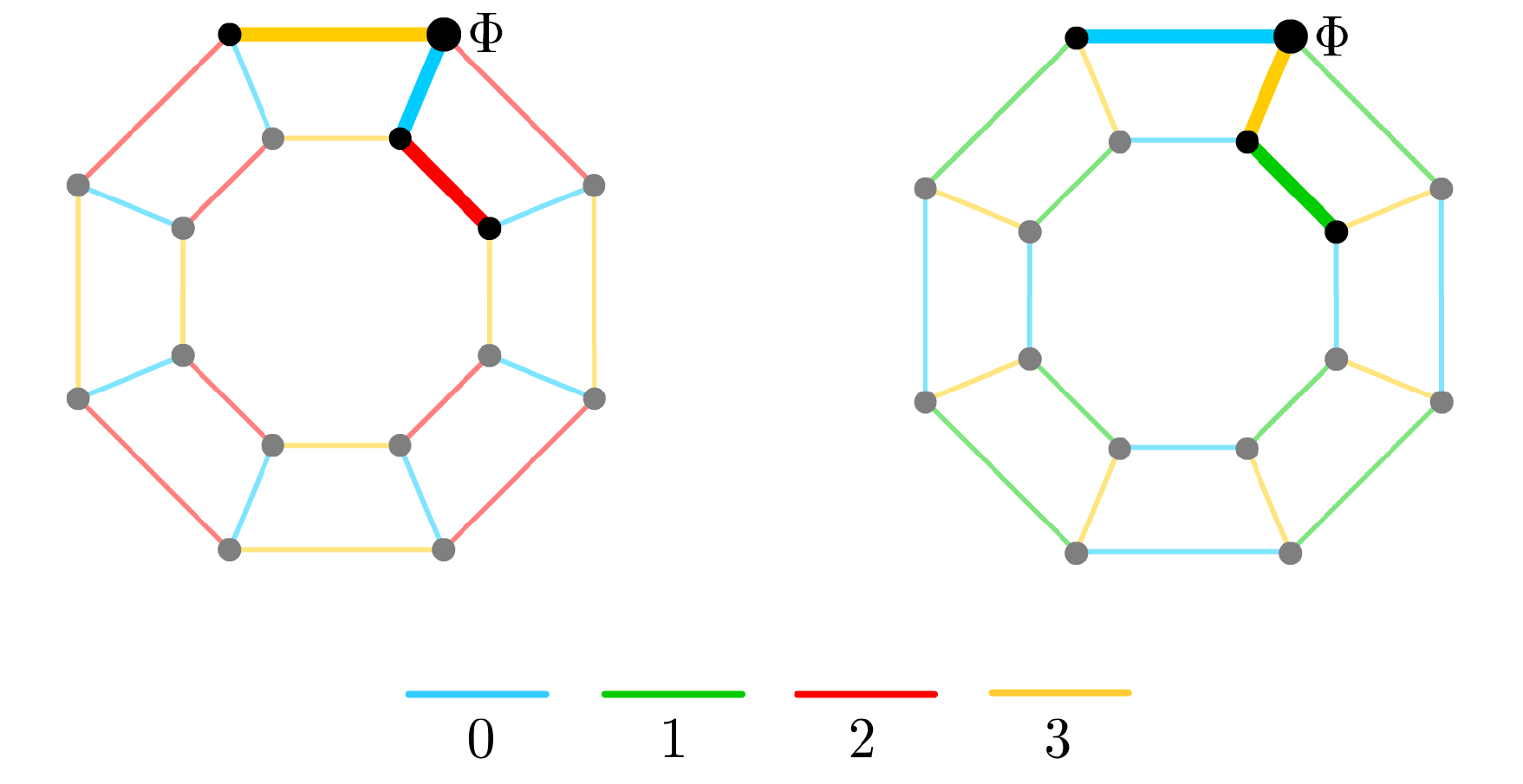}
    \caption{A $1$-face and a $2$-face of $\B$, with their respective intersections with $E_\Theta$ shown in thick edges. All $1$- and $2$-faces of $\B$ intersect $E_{\Theta}$ in a similar way.}
    \label{fig:2-face}
\end{figure}

\begin{figure}[H]
\begin{center}
\includegraphics[width=1\textwidth]{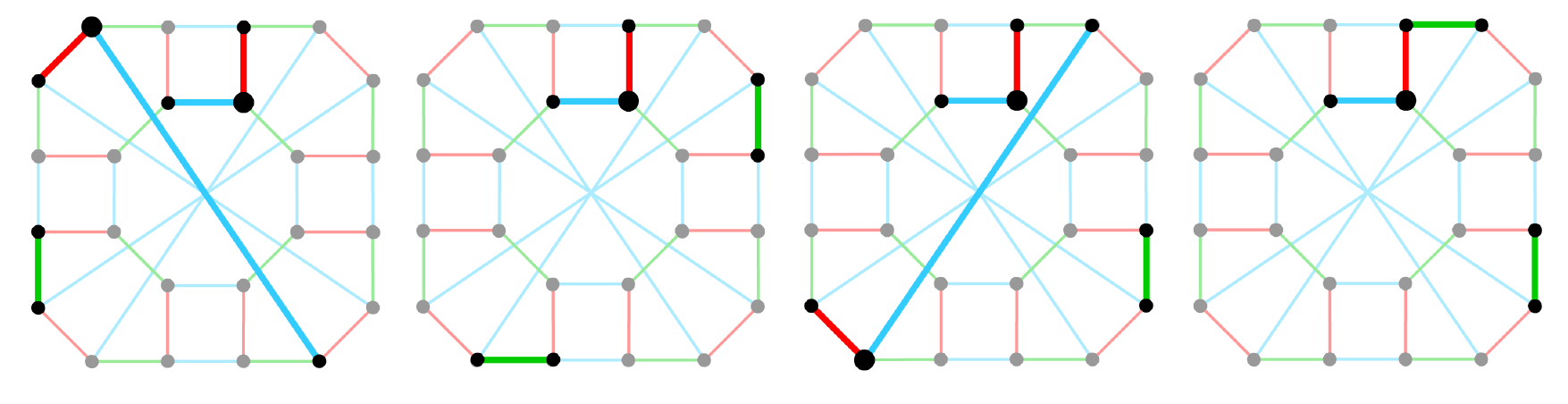}
\caption{The four facets of $\B$, with their intersection with $E_\Theta$ shown in thick edges.}
\label{fig:facets}
\end{center}
\end{figure}

As we will see in the following pages, the cover $\Cov(\B,\zeta_{\Theta})$ arising from this voltage assignment is an unfaithful polytopal maniplex. 
In particular, we show in Lemmas \ref{lem:12face} and \ref{lem:03face} that if $F$ is an $i$-face of $\B$ for some $i\in \{0,\ldots,3\}$, then $\pi^{-1}(F)$ is an $i$-face of $\Cov(\B,\zeta_{\Theta})$. 
For this, it suffices to show that the set edges of $E_{\Theta}$ contained in $F$ do not form a cut-set of $F$ by Lemma \ref{lem:voltcon}. 
The visually inclined reader can convince themselves that this holds by inspecting Figures \ref{fig:2-face} and \ref{fig:facets}, showing how the set $E_\Theta$ intersects each $1$-, $2$ and $3$-face of $\B$, respectively. 
The self-duality of $\B$ guarantees that the same holds for each $0$-face. 
Proofs of these two lemmas are included for the sake of completion.
However, those readers who are convinced by the figures may skip (the auxiliary) Lemma \ref{lem:bes} and the proofs of Lemmas \ref{lem:12face} and \ref{lem:03face}.

\begin{lemma}
\label{lem:bes}
Let $F$ be an $i$-face of $\B$. If $i \in \{1,2\}$ then all the following hold:
\begin{enumerate}
\item[(B.1)] If the $i$-edge between $\Phi$ and $\Psi$ is in $E_\Theta$, then either $\Phi \in \Theta^{i+2}$ or $\Psi \in \Theta^{i+2}$, where addition is taken modulo $4$;
\item[(B.2)] If the $i$-edge between $\Phi$ and $\Psi$ is in $E_\Theta$, then $\{\Phi,\Psi\} \cap \Theta = \emptyset$.
\item[(B.3)] There is exactly one $j$-edge in $E_\Theta \cap F$ for each $j \in \{0,1,2,3\}\setminus\{i\}$.
\end{enumerate}

Moreover if $i \in \{0,3\}$ then one of the following holds: 
\begin{enumerate}
\item[(B.4)] $F$ contains exactly two $j$-edges in $E_\Theta$ for each $j \in \{i+1,i-1\}$ and one $(i+2)$-edge in $E_\Theta$;
\item[(B.4)'] $F$ contains exactly one $j$-edge in $E_\Theta$ for each $j \in \{i+1,i-1\}$ and two $(i+2)$-edges in $E_\Theta$;
\end{enumerate}
\end{lemma}

\begin{proof}
Items (B.1) and (B.2) follow from the definition of $E_\Theta$.

For item (B.3), let us assume that $i=2$ (a similar argument applies for $i=1$).
First, let $\Phi$ be the flag in $\Phi\in F\cap \Theta$; it exists and is unique by (A.1). 
Since all the $0$- and $3$-edges of $E_\Theta$ have one flag in $\Theta$, the only $0$- and $3$-edges in $E_\Theta\cap F$ are precisely $(\Phi,\Phi^0)$ and $(\Phi,\Phi^3)$.
Now, the $1$-edges of $E_\Theta$ all have the form $(\Psi^3,\Psi^{31})$, for some $\Psi\in \Theta$. 
Note that $\Psi$ and $\Psi^3$ are in the same $2$-face, so if $(\Psi^3,\Psi^{31})$ is a $1$-edge in $E_\Theta\cap F$, then $\Psi\in \Theta\cap F$, and therefore again by (A.1) $\Psi=\Phi$, implying that there is only $1$-edge in $E_\Theta\cap F$.

Finally, observe that if $i \in \{0,3\}$ then one of the (A.3) and (A.4) holds for $F$. 
It is straightforward to see that if (A.3) holds, then so does (B.4).
Similarly, if (A.4) holds, then (B,4)' holds as well.
\end{proof}

\begin{lemma}
\label{lem:12face}
If $F$ is an $i$-face of $\B$ with $i \in \{1,2\}$, then $\pi^{-1}(F)$ is a connected subgraph of $\Cov(\B,\zeta_\Theta)$.
\end{lemma}

\begin{proof}
By (B.3) $F$ contains exactly one $j$-edge of $E_\Theta$ for each $j \in \{0,1,2,3\}\setminus\{i\}$. 
As $F$ contains no $i$-edges, those three edges are the only edges in $F \cap E_\Theta$.
Let $F^-$ be the graph obtained by removing the edges in $F\cap E_\Theta$  from $F$.
Note that $F^-$ remains connected (see \Cref{fig:2-face}).
Since every edge lying on $F^-$ has trivial voltage, by \Cref{lem:voltcon}  $\pi^{-1}(F)$ is connected. 
\end{proof}

\begin{lemma}
\label{lem:03face}
If $F$ is an $i$-face of $\B$, $i \in \{0,3\}$, then $\pi^{-1}(F)$ is a connected subgraph of $\Cov(\B,\zeta_\Theta)$.
\end{lemma}

\begin{proof}
There are two cases to consider, depending on which of the conditions (B.4) and (B.4)' is satisfied by $F$. For the sake of simplicity, we may assume $i = 3$, and the case $i = 0$ follows by duality.

First, suppose that $(B.4)$ holds (and thus also (A.3)). 
That is, there are exactly two flags $\Phi,\Psi \in F \cap X$.
Furthermore, there are exactly five edges in $F \cap E_\Theta$: the $0$- and $2$-edges incident to $\Phi$ and $\Psi$ and a single $1$-edge.
Let $F'$ be the graph resulting from the deletion of the $0$- and $2$-edges in $F \cap E_\Theta$. 
Note that $F'$ can only be disconnected if $\Phi$ and $\Psi$ are the endpoints of a $1$-edge. 
However, this cannot happen because then $\Phi$ and $\Psi$ would lie on the same $2$-face, contradicting (A.3). 
Now, the deletion of a $1$-edge $e$ disconnects $F'$ only if one of $\Phi$ or $\Psi$ is an endpoint of $e$, but this cannot happen either, since, by (B.2), no endpoint of $e$ is an element of $\Theta$.
It follows that the graph obtained from $F$ after the deletion of the edges in $F \cap E_\Theta$ is connected. Therefore, $\pi^{-1}(F)$ is connected by \Cref{lem:voltcon}. 

Now, suppose $(B.4')$ holds. Then $F \cap E_\Theta$ consists of two $1$-edges and one $j$-edge for each $j \in \{0,2\}$. 
Moreover, the $0$- and $2$-edges share a common endpoint, $\Phi$, which is an element of $\Theta^0$. 
Let $F'$ be the graph resulting from the removal of the edges in $F \cap E_\Theta$ from $F$. 
Again, $F'$ can only be disconnected if $\Phi$ is an endpoint of one of the $1$-edges in $F \cap E_\Theta$, but according to the definition of $\zeta_\Theta$, no flag of $\B$ is incident to three edges in $E_\Theta$.
\end{proof}


\begin{lemma}
\label{lem:Bstarmani}
The covering graph $\B^*:=\Cov(\B,\zeta_\Theta)$ is a maniplex.
\end{lemma}

\begin{proof}
We need to prove that $\B^*$ is connected and that every $ijij$-path is a cycle of length $4$ whenever $|i-j| >1$. Let us begin by proving that $\B^*$ is connected.
Recall that the flags of $\B^*$ are of the form $(\Phi,r)$ where $\Phi$ is a flag of $\B$ and $r \in \{0,1\}$. 

By Lemmas \ref{lem:12face} and \ref{lem:03face}, we see that if $(\Phi,r)$ is a flag of $\B^*$, then $\B^*_{\bar{i}}((\Phi,r)) = \pi^{-1}(\B_{\bar{i}}(\Phi))$. 
Moreover, two subgraphs $\B^*_{\bar{i}}((\Phi,r))$ and $\B^*_{\bar{j}}((\Phi',r'))$ have non-trivial intersection if and only if the corresponding projections $\B_{\bar{i}}(\Phi)$ and $\B_{\bar{j}}(\Phi')$ have non-trivial intersection.
The connectedness of $\B^*$ then follows from the fact that $\B^*_{\bar{i}}((\Phi,r)) = \pi^{-1}(\B_{\bar{i}}(\Phi))$ for all $(\Phi,r) \in \B^*$ and $i\in\{0,1,2,3\}$ and that $\B$ is connected.

Now, let $C$ be an $ijij$-path in $\B$ with $|i-j|>1$.  
Then, there are only two options for the voltages of the edges in $C$: either all edges in $C$ have trivial voltage, or two consecutive edges (and only these two) have trivial voltage. 
In either case, the number of edges with non-trivial voltage is even and by Lemma \ref{lem:voltmani}, we conclude that $\B^*$ is a maniplex.
\end{proof}

Lemmas \ref{lem:12face}, \ref{lem:03face} and \ref{lem:Bstarmani} combined give us:

\begin{proposition}\label{prop:sameposet}
The covering graph $\Cov(\B,\zeta_\Theta)$ is a $4$-maniplex with the same poset as $\B$.
\end{proposition}

We have just shown that $\B^*=\Cov(\B,\zeta_\Theta)$ is a maniplex that has flag-set $\B \times \mathbb{Z}_2$. That is, $\B^*$ has twice as many flags as $\B$.
\Cref{prop:sameposet} implies that $\Pos(\B^*)=\Pos(\B)$, which in turn implies that the number of flags of $\Pos(\B^*)$ is half the number of flags of $\B^*$ and therefore $\B^*$ is unfaithful.
Moreover, $\Pos(\B^*)=\Pos(\B)=\po$, which means that $\B^*$ is a polytopal maniplex.
Therefore, we have established the following results.

\begin{proposition}
    $\B^*$ is an unfaithful polytopal $4$-maniplex.
\end{proposition}

The above result and \Cref{teo:theequivalnce} imply that:

\begin{theorem}\label{teo:rank4}
    In rank $4$ there exists at least one sparse group that is not semisparse.
\end{theorem}

\section{Extensions}\label{sec:extensions}
We have shown that there exists a $4$-maniplex that is unfaithful and polytopal. 
In this section we study an extension of maniplexes that given an unfaithful and polytopal $n$-maniplex constructs an $(n+1)$-maniplex that is also unfaithful and polytopal, thus showing that for all $n>3$, there exist sparse groups that are not semisparse.

Given an $n$-maniplex $\M$ there are two classical ways to construct an $(n+1)$-maniplex whose facets are all isomorphic to $\M$: the trivial extension and the $2^{\widehat{\M}}$ (see \cite{douglas2018twist}). 
These are, in fact, two extensions that can be seen as a family of {\em colour-coded extensions}, described in \cite[Section 6.1]{douglas2018twist}. 
Here, we analyse one particular instance of these colour-coded extensions.

Let $\M$ be an $n$-maniplex with monodromy group $\Mon(\M) = \langle r_0, \ldots, r_{n-1} \rangle$. Let $F$ be a fixed facet of $\M$.
Define the $(n+1)$-maniplex $\M^F$ as follows. The set of flags of $\M^F$ is the set of pairs $(\Phi,x)$ where $\Phi \in \F(\M)$ and $x \in \ZZ_2^2$. 
The $i$-edges of $\M^F$ are the edges between a flag $(\Phi,x)$ and the flag $s_i(\Phi,x)$, where $s_0, s_1, \dots, s_{n}$ are the generators of the monodromy group $\Mon(\M^F)$, and are defined as:

\[
s_i(\Phi,x)=
 \begin{cases}
(r_i\Phi,x) & \text{ if $0 \leq i \leq n-1$}, \\ 
(\Phi,x + (1,0)) &\text{ if $i = n$ and $\Phi \notin F$}, \\
(\Phi,x + (1,1)) &\text{ if $i = n$ and $\Phi \in F$}.
\end{cases}
\]

Then $\M^F$ is isomorphic to the colour-coded extension $2^{(\M,C)}$ where $C := \{c_1,c_2\}$ is a set of two colours and the facet $F$ has colour $c_1$ while all the other facets have colour $c_2$.

Given a flag $(\Phi,x)$ of $\M^F$, the facet $\M^F_{\bar{n}}(\Phi,x)$ consists of all the flags $w(\Phi,x)$, were $w\in \langle s_0, s_1, \dots, s_{n-1} \rangle$, together with all the $i$-edges joining such flags, with $i\neq n$.
Note that by definition, the group $ \langle s_0, s_1, \dots, s_{n-1} \rangle$ is transitive on the flags of $\M$, and no element of such a group changes the second coordinate,  implying that $\M^F_{\bar{n}}(\Phi,x) = \{(\Psi, x) \mid \Psi \in \M \}$.
Therefore, $\M^F$ has exactly four facets, one corresponding to each of the elements of $\ZZ_2^2$.

\begin{proposition}
\label{prop:extunfaithful}
If $\M$ is unfaithful, then so is $\M^F$.
\end{proposition}

\begin{proof}
Let $\Phi_1$ and $\Phi_2$ be two flags of $\M$ such that there exist $i$-faces $F_i$ with $\Phi_1, \Phi_2 \in \bigcap_{i=0}^{n-1} F_i$. 
For each $i \in \{0,\ldots,n\}$, let $G_i$ be the $i$-face of $\M^F$ containing the flag $(\Phi_1,(0,0))$. That is, $G_i = \M^F_{\bar{i}}(\Phi_1,(0,0))$. 

Given $i \in \{0,\ldots,n-1\}$, since $\Phi_1, \Phi_2 \in F_i$, there exists $w_i \in \langle r_0,\ldots,r_{i-1},r_{i+1},\ldots,r_{n-1} \rangle$ such that $w_i\Phi_1 = \Phi_2$. 
In other words, there is a path from $\Phi_1$ to $\Phi_2$ using only colours in $\{0,\ldots,i-1,i+1,\ldots,n-1\}$. 
Suppose $w_i = r_{i_0}r_{i_1}\ldots r_{i_{k(i)}}$ for some integer $k(i)$ (that depends on $i$) and let $\bar{w}_i = s_{i_0}s_{i_0}\ldots,s_{i_{k(i)}}$ be the corresponding element in $\Mon(\M^F)$. 
Then, since $w_i$ never changes the second coordinate, $(\Phi_2,(0,0)) = (w_i\Phi_1,(0,0)) = \bar{w}_i(\Phi_1,(0,0)) \in G_i$.
It follows that both $(\Phi_1,(0,0))$ and $(\Phi_2,(0,0))$ are in $G_i$, for each $i \in \{0,\ldots,n-1\}$, and thus they are both in $\bigcap_{i=0}^{n-1} G_i$.

Moreover, none of the paths from $(\Phi_1,(0,0))$ to $(\Phi_2,(0,0))$ described above trace edges of colour $n$, implying that $(\Phi_2,(0,0)) \in G_n$. Therefore, both $(\Phi_1,(0,0))$ and $(\Phi_2,(0,0))$ are in $\bigcap_{i-0}^n G_i$. We conclude that $\M^F$ is unfaithful.
\end{proof}

Our goal in this section is to show that the extension $\M^F$ is polytopal whenever $\M$ is polytopal. 
For this, it will be useful to have a description of the faces of $\M^F$ that intersect a given facet $G_n$. 
First note that if $G_n$ is a facet of $\M^F$ and contains the flag $(\Phi, x)$, then there is a path with colours in $\{0,1, \dots, n-1\}$ from any other flag of $G_n$ to $(\Phi,x)$.
That is, if $(\Psi, y)$ is another flag of $G_n$, then there exists $w\in \langle s_0, s_1, \dots, s_{n-1} \rangle$ such that $w(\Phi, x)=(\Psi, y)$. 
But by the definition of $s_i$, this implies that $x=y$, and thus $G_n \cong \M$. 
Note that this does not imply that the section $G_n/G_{-1}$ of $\Pos(\M^F)$ is isomorphic to the poset $\Pos(\M)$.
This is because the elements of $G_n/G_{-1}$ are subgraphs of $\M^F$, and thus they have edges of colour $n$, which the subgraphs of $\M$ do not have.
We will see that, in fact, whenever $\M$ has the diamond condition, $G_n/G_{-1}$ and $\Pos(\M)$ are indeed isomorphic posets.
%
Without loss of generality, we may assume that $G_n = \{(\Psi,(0,0)) \mid \Psi \in \M \}$. 
As pointed out before, we want to describe the $i$-faces of $\M^F$ that have non-empty intersection with $G_n$.

Let $i <n$, let $\Phi \in \M$, and define
\[
Y_{i,\Phi} = 
\begin{cases}
\{(0,0),(1,0)\} &\text{ if } \M_{\bar{i}}(\Phi) \cap F = \emptyset, \\
\{(0,0),(1,1)\} &\text{ if } \M_{\bar{i}}(\Phi) = F, \\
\ZZ_2^2 &\text{ if } \M_{\bar{i}}(\Phi) \cap F \neq \emptyset \text{ but } \M_{\bar{j}}(\Phi) \not \subset F.  \\
\end{cases}
\]

Observe that if $\M$ satisfies the diamond condition, then no $i$-face of $\M$, $i < n-1$, can be properly contained in $F$, as it would be incident to only one facet. Therefore, if $\M$ is polytopal, then the three possibilities for $\M_{\bar{i}}(\Phi)$ in the definition of $Y_{i,\Phi}$ above are exhaustive and $Y_{i,\Phi}$ is well defined. 
Now, let 
\begin{align*}
V_{i,\Phi} = \{ (\Psi, y) \mid \Psi \in \M_{\bar{i}}(\Phi), y \in Y_{i,\Phi}\},
\end{align*}
and let $G_i(\Phi)$ be the subgraph of $\M^F_{\bar{i}}$ induced by the set of flags $V_{i,\Phi}$. 
Observe that $G_i(\Phi)$ is connected due to the definitions of $r_n$ and $Y_{i,\Phi}$, and the fact that $\M_{\bar{i}}(\Phi)$ is connected. 
Moreover, if for some $j\neq i$ there exists a flag $(\Psi,x)$ that is $j$-adjacent to a flag $(\Phi',y) \in G_i(\Phi)$, then the flags $\Psi$ and $\Phi'$ are $j$-adjacent in $\M$ or $j=n$. 
In both cases, we have $\M_{\bar{i}}(\Psi) = \M_{\bar{i}}(\Phi')$. It follows that $Y_{i,\Psi} = Y_{i,\Phi'}$  and thus $(\Psi,x) \in G_i(\Phi)$. 
That is, $G_i(\Phi)$ is in fact an $i$-face in $G_n/G_{-1}$, and therefore is an $i$-face of $\M^F$. 
Conversely, for every $i$-face $G_i'$ in $G_n/G_{-1}$, there exists a flag $\Phi \in \M$ for which $G_i' = G_i(\Phi)$. 

With this description of the faces in $G_n/G_{-1}$, it is not difficult to show that $G_n/G_{-1}$ is, in fact, isomorphic to $\Pos(\M)$. 
Let $\varphi \colon \Pos(\M) \to G_n/G_{-1}$ be the function mapping each $i$-face $\M_{\bar{i}(\Phi)}$ to $G_i(\Phi)$. We claim that $\varphi$ is an isomorphism. 

 We first need to show that $\varphi$ is a well-defined function. Suppose $\M_{\bar{j}}(\Phi) = \M_{\bar{k}}(\Phi')$ for some $j,k < n$ and some $\Phi,\Phi' \in \M$. 
 Then $j = k$  and since $\M_{\bar{j}}(\Phi) = \M_{\bar{j}}(\Phi')$, the definition of $G_j(\Phi)$ implies that $G_j(\Phi) = G_k(\Phi')$.
 Then $\varphi$ does not depend on the choice of flag representatives and, thus, is well defined.
 Now, $\varphi$ is clearly surjective, as every $i$-face of $G_n/G_{-1}$ is of the form $G_i(\Phi)$ for some $\Phi \in \M$ and $G_i(\Phi) = \varphi(\M_{\bar{i}}(\Phi))$.
 To see that $\varphi$ is injective, let $\M_{\bar{i}}(\Phi)$ and $\M_{\bar{j}}(\Phi')$ be two faces of $\M$ and suppose $G_i(\Phi) = G_j(\Phi')$. 
 Then $i = j$, and the flags $(\Phi,(0,0))\in G_i(\Phi)$ and $(\Phi',(0,0))\in G_i(\Phi')$ are in the same $i$-face of $\M^F$, implying that $\Phi$ and $\Phi'$ are in the same $i$-face of $\M$. 
Therefore, $\M_{\bar{i}}(\Phi) = \M_{\bar{i}}(\Phi') = \M_{\bar{j}}(\Phi')$ from which we see that $\varphi$ is injective.
 
 Finally, since incidence is given by intersections, $\varphi$ maps incident faces to incident faces and non-incident faces to non-incident faces. We have just proved the following.
 
\begin{lemma}
\label{lem:seccioniso}
If $\M$ satisfies the diamond condition and $G_n$ is a facet of $\M^F$, then $G_n/G_{-1}$ is isomorphic to $\Pos(\M)$. 
\end{lemma}

The above lemma has the following two important corollaries, which are instrumental in proving that $\M^F$ is a polytopal maniplex.

\begin{corollary}
\label{cor:extdiamante}
Let $\M$ be a polytopal maniplex of rank $n$, and let $G_n$ be a facet of $\M$.
Then, the section $G_n/G_{-1}$ of the poset $\Pos(\M^F)$ satisfies the diamond condition.
\end{corollary}

\begin{corollary}
\label{cor:extconexidad}
Let $\M$ be a polytopal maniplex of rank $n$, let $i < n$ and let $G_i$ and $G_n$ be incident faces of $\M^F$ of ranks $i$ and $n$, respectively. Then, the section $G_n/G_i$ is flag connected.
\end{corollary}

Corollary \ref{cor:extconexidad} will be useful later in proving that $\M^F$ is strongly flag connected. For the moment, observe that by Corollary \ref{cor:extdiamante} we only need to show that every face of rank $n-1$ of $\Pos(\M^F)$ is incident to exactly two faces of rank $n$ to show that $\Pos(\M^F)$ satisfies the diamond condition.

\begin{lemma}\label{lemma:n-diamond}
Each $(n-1)$-face of $\Pos(\M^F)$ is incident to exactly two $n$-faces of $\Pos(\M^F)$.
\end{lemma}

\begin{proof}
Let $G_{n-1}$ be an $(n-1)$-face of $\Pos(\M^F)$, let $\tilde{\Phi} := (\Phi,x) \in G_{n-1}$ and let $F_{n-1} = \M_{\overline{n-1}}(\Phi)$. 

We have two cases depending on whether $F_{n-1}$ intersects $F$ or not.
Recall that $F$ is an $(n-1)$-face of $\M$, and thus whenever $F_{n-1}\cap F \neq \emptyset$, we have in fact $F = F_{n-1}$. 
We then see that $G_{n-1}$ has more flags than $F$, since it has $n$-edges but $F$ does not.

So suppose $F_{n-1}\cap F \neq \emptyset$, and hence $F_{n-1}=F$; in particular this implies that $\Phi \in F$, and therefore the flags of $G_{n-1}(\Phi)$ are those in the set $V_{n-1,\Phi}=\{(\Psi, y) \mid \Psi\in F, y \in \{(0,0),(1,1)\}$. 
Then, a path starting at $(\Phi,x)$ and tracing only edges of colours in $\{0,\ldots, n-2\}$ must end on a flag of the form $(\Psi,x)$ with $\Psi \in F$. 
What is more, the $n$-neighbour of such a flag is $(\Psi,x+(1,1))$, implying that $G_{n-1}$ only has elements whose second coordinate is $x$ or $x + (1,1)$. 
In other words, $G_{n-1}$ is incident only to the two $n$-faces of $\M^F$ having second coordinate either $x$ or $x+(1,1)$.

Now suppose $F_{n-1}\cap F = \emptyset$.
Then $V_{n-1,\Phi} = \{(\Psi,y) \mid \Psi \in F_{n-1}, y\in\{(0,0),(1,0)\}\}$.
Again, a path starting at $(\Phi,x)$ and tracing only edges of colours in $\{0,\dots,n-2\}$ must end in a flag of the form $(\Psi,x)$ with $\Psi\in F_{n-1}$. 
Thus, a path with colour in $\{0,\ldots,n-2,n\}$ starting at $(\Phi,x)$ ends in a flag of the form $(\Psi,y)$ with $\Psi\in F_{n-1}$, and $y=x$ or $y = x+(0,1)$.
Therefore, the two $n$-faces incident to $G_{n-1}$ are precisely the two $n$-faces with flags having second coordinate $x$ and $x+(0,1)$, respectively.
\end{proof}




\begin{corollary}
\label{cor:extdiam2}
If $\Pos(\M)$ satisfies the diamond condition, then so does $\Pos(\M^F)$.
\end{corollary}

We are now ready to prove that $\Pos(\M^F)$ is strongly flag-connected.
That is, we need to show that each section $\Pos(\M^F)$ is flag-connected. 
Note that Corollary \ref{cor:extconexidad} does most of the heavy lifting here.
Indeed, by \Cref{cor:extconexidad} we need only to show that for each $i$-face $G_i$ with $i \in \{-1, \ldots, n-1\}$, the section $G_{n+1}/G_i$ is flag connected, where $G_{n+1}$ is the unique maximal face of $\M^F$.
Furthermore, by \Cref{lemma:n-diamond}, the sections $G_{n+1}/G_{n-1}$ are also flag connected (they have exactly two flags), so we are left with the case where $i \in \{-1, \ldots, n-2\}$ to be dealt with.

\begin{lemma}
Let $i \in \{-1, \ldots, n-2\}$ and let $G_{n+1}$ be the unique maximal face of $\M^F$. 
Then, the section $G_{n+1}/G_i$ is flag-connected.
\end{lemma}

\begin{proof}
Let $C$ and $C'$ be two maximal chains contained in $G_{n+1}/G_i$. We must show that there is a sequence of maximal chains in $G_{n+1}/G_i$ connecting $C$ to $C'$ and such that consecutive elements are adjacent chains. Let $G_n$ and $G_n'$ be the elements of rank $n$ of $C$ and $C'$, respectively. If $G_n = G_n'$, then by Corollary \ref{cor:extconexidad}, we are done. 
Therefore, we assume that $G_n \neq G_n'$ and suppose, as we may, that $G_n = \{(\Psi,(0,0)) \mid \Psi \in \M\}$, and that, for some $x \in \ZZ^2_2 \setminus \{(0,0)\}$, $G_n'= \{(\Psi,x) \mid \Psi \in \M\}$.

Suppose that there exists a flag $(\Phi, (0,0))\in G_n \cap G_i$ such that $s_n(\Phi, (0,0)) \in G_n'$. 
By the definition of $s_n$, this is equivalent to having that $x=(1,0)$ and $\Phi \notin F$, or $x=(1,1)$ and $\Phi \in F$.
In this case, we will construct maximal chains $C_\Phi$ and $C'_\Phi$ that are $n$-adjacent and satisfying that $G_n \in C_\Phi$, and $G_n'\in C'_\Phi$; these chains will give us sequences of adjacent flags that we then concatenate to obtain the desired sequence from $C$ to $C'$.

Observe that for all $i \leq j<n$, the $j$-face $\M^F_{\bar{j}}(\Phi,(0,0))$ is in $G_{n}/G_i$ while $\M^F_{\bar{j}}(\Phi,x)$ is in $G_{n}'/G_i$. 
Furthermore, since $(\Phi,(0,0))$ is $n$-adjacent to $(\Phi,x)$ we have $ \M^F_{\bar{j}}(\Phi,(0,0))= \M^F_{\bar{j}}(\Phi,x)$. For $i \leq j<n$ let $G_j$ denote the $j$-face $\M^F_{\bar{j}}(\Phi,(0,0))$. 
Then the chains $C_\Phi:=\{G_i,G_{i+1},G_{i+2} ,\ldots,$ $G_{n-1},G_n\}$ and $C_\Phi':=\{G_i,G_{i+1},G_{i+2} ,\ldots,G_{n-1},G_n'\}$ are maximal chains of $G_n/G_i$ and $G_n'/G_i$, respectively, and differ only in their elements of rank $n$. 

Since $G_n/G_i$ is flag connected (by Corollary \ref{cor:extconexidad}), there exists a sequence $S$ of maximal chains of $G_n/G_i$ connecting $C \setminus \{G_{n+1}\}$ to $C_{\Phi}$ and such that consecutive elements are adjacent. Similarly, there is a sequence $S'$ of maximal chains of $G_n'/G_i$ connecting $C_{\Phi}'$ to $C'\setminus \{G'_{n+1}\}$ and such that consecutive elements are adjacent. Then, the concatenation of the sequences $S$ and $S'$ is a sequence of maximal chains in $G_{n+1}/G_i$ (after adding the unique element $G_{n+1}$ to each) connecting $C$ and $C'$ and such that consecutive elements are adjacent.

Now we are only missing the case when $x=(0,1)$. In that case we consider the $n$-face $G_n''= \{(\Psi,(0,1)) \mid \Psi \in \M\}$, and note that for $\Psi \in F$ and $\Psi' \notin F$, then $s_n(\Psi, (0,1))=(\Psi, (1,0)) \in G_n$, and $s_n(\Psi', (0,1))=(\Psi', (1,1)) \in G_n'$. 
And by the previous case, one can construct sequences of adjacent maximal chains and then concatenate them to obtain the sequence between $C$ and $C'$.
This concludes the proof.  
\end{proof}

\begin{corollary}
\label{cor:extcon2}
If $\M$ is a polytopal maniplex of rank $n$, then $\M^F$ is strongly flag connected.
\end{corollary}

The combination of Proposition \ref{prop:extunfaithful} with Corollaries \ref{cor:extdiam2} and \ref{cor:extcon2} gives us the following proposition.

\begin{proposition}
\label{prop:extentiongen}
If $\M$ is an unfaithful polytopal maniplex of rank $n$, then $\M^F$ is an unfaithful polytopal maniplex of rank $n+1$.
\end{proposition}

We are now ready to state the main result of this paper, which disproves Hartley's conjecture \cite[Conjecture 5.2]{hartley1999more}.
The theorem is a result of the above proposition, \Cref{teo:theequivalnce} and \Cref{teo:rank4}.

\begin{theorem}
    For all $n>3$, there exists a group $N\leq \C$ that is sparse but not semisparse.
\end{theorem}

\section{Concluding remarks}\label{sec:remarks}
We have constructed an infinite family of unfaithful polytopal maniplexes having one element of every rank $n>3$. We have done so by first constructing one such maniplex of rank $4$ (namely, $\B^*$) and repeatedly applying the extension construction of \Cref{sec:extensions} to obtain one maniplex for every rank greater than $4$. Since in rank $3$ every sparse group is semisparse, the family we constructed is, in a sense, as good as possible when it comes to showing that sparse groups that are not semisparse are abundant and exist `for every possible rank'. However, this is not the only family of unfaithful polytopal maniplexes that can be constructed using the techniques in this paper.

Consider for instance the regular abstract polytope $\po$ of rank $6$ denoted by $\{4,3,6,3,4\}^*1728$ in Hartley's atlas of small regular polytopes \cite{hartley2006atlas}. As it transpires, $\po$ has the property that for every $i \in \{0,\ldots,5\}$, every $i$-face is non-orientable (equivalently, every $i$-face of the flag-graph of $\po$ is non-bipartite). We may call this property strong non-orientability. Observe that if $\M$ is the flag-graph of a strongly non-orientable polytope, then its canonical double cover $\overline{\M}$ is an unfaithful polytopal maniplex. Indeed, every $i$-face of $\M$ lifts to a connected graph and is thus an $i$-face of $\overline{\M}$, from which we see that both $\M$ and $\overline{\M}$ have isomorphic posets even though $\overline{\M}$ is twice as large as $\M$. We can then construct extension of $\overline{\M}$ of higher ranks by using \Cref{prop:extentiongen}. That is, we can obtain an infinite family of unfaithful polytopal maniplexes from every strongly non-orientable polytope. In particular, $\po$ yields such a family. We believe, even though we make no claim, that strongly non-orientable polytopes are common for ranks $n \geq 6$ (no such polytope exists for ranks smaller than $6$).

Finally, we would like to point out that we considered double covers in Section \ref{sec:bigconstruction} merely because they are easy to handle, but in principle any $n$-fold cover of (the flag-graph of) a polytope where every $i$-face lifts to an $i$-face produces an unfaithful polytopal maniplex, and by the results in \Cref{sec:extensions}, an entire infinite family.

\section*{Acknowledgements}
The first author gratefully acknowledges the financial support of  of CONACyT grant A1-S-21678 and PAPIIT- DGAPA grant IN109023. 

The second author gratefully acknowledges financial support from the F\'ed\'eration Wallonie-Bruxelles -- Actions de Recherche Concert\'ees (ARC Advanced grant).

\bibliographystyle{amsplain}
\bibliography{sparse}

\end{document}